\documentclass[leqno,10pt]{article}
\usepackage{epsfig}
\usepackage{graphicx}
\usepackage{CJKutf8}
\usepackage{amssymb}
\usepackage{amsmath}
\usepackage{authblk}
\usepackage{setspace}

\def\qed{\vrule height 1.2ex width 1.2ex depth 0ex }
\newcommand{\dis}{\displaystyle}
\newcounter{corollary}
\newtheorem{refproof}{Proof of corollary}
\newtheorem{theorem}{\bf Theorem}
\newtheorem{Lemma}{\bf Lemma}
\newtheorem{Proposition}{\bf Proposition}
\newtheorem{corollary}{Corollary}

\newtheorem{definition}{\bf Definition}
\usepackage[colorlinks=true]{hyperref}
\hypersetup{urlcolor=blue, citecolor=red}
\allowdisplaybreaks
\usepackage{geometry}\usepackage[super]{nth}

\newenvironment{corollaryenv}{\begin{corollary}$\ $}{\end{corollary}\vspace{\baselineskip}}

\title{}
\author{}
\date{}
\begin{document}
\maketitle
%%%%%%%%%%%%%%%%%%%%%%%%%%%%%%%%%
\thispagestyle{plain}
\begin{center}
    \begin{huge}
    \setstretch{1}{\bf{On the second order of Zeta functional equations for Riemann Type}}\\
    \end{huge}
\end{center}
\vspace{6pt}
%%%%%%%%%%%%%%%%%%%%%%%%%%%%%%%%%
\begin{flushleft}
\begin{large}
\bf{Chin-yuan Hu $^1$, Tsung-lin Cheng $^{2,*}$ and Ie-bin Lian $^{3,*}$}
\end{large}
\end{flushleft}
\vspace{6pt}
%%%%%%%%%%%%%%%%%%%%%%%%%%%%%%%%%
\begin{flushleft}
$^1$ National Changhua University of Education, Changhua, 50002, Taiwan\\
\emph{E-mail: cyhcyh1101@icloud.com}
\end{flushleft}
%%%%%%%%%%%%%%%%%%%%%%%%%%%%%%%%%
\begin{flushleft}
$^2$ Graduate Institute of Statistics and Information Science, National Changhua\\ University of Education, Changhua, 50002,Taiwan\\
\emph{E-mail: tlcheng@cc.ncue.edu.tw}
\end{flushleft}
%%%%%%%%%%%%%%%%%%%%%%%%%%%%%%%%%
\begin{flushleft}
$^3$ Graduate Institute of Statistics and Information Science, National Changhua\\
University of Education, Changhua, 50002,Taiwan\\
\emph{E-mail*: maiblian@cc.ncue.edu.tw}
\end{flushleft}
%%%%%%%%%%%%%%%%%%%%%%%%%%%%%%%%%
AMS subject classification: 11M06, 60J65, 60E07\\
Keywords: Riemann Zeta function, Mellin transformation characteristic function,
functional equation.
%%%%%%%%%%%%%%%%%%%%%%%%%%%%%%%%%%%%%%%%%%%%%%%%%%%%%%%%%%%%%%%%%%%%%%%%%%
\newpage
\newgeometry{%版面設定
total={180mm,270mm}, 
top=40mm,
bottom=40mm,
left=30mm,
right=30mm}

\begin{abs}
\begin{normalsize}
This paper discuss a new class of functional equations by using both Poisson’s
summation formula and Jacobi’s type theta function. The class of Riemann’s type
functional equations are derived from self-reciprocal probability density functions.
Finally, the second order Zeta functional equations for Riemann’s type is also
investigated.
\end{normalsize}
\end{abs}
\section{Introduction}
This paper is concerned with Riemann’s Zeta function $\zeta(s)$, which is defined by
\begin{equation*}
    \zeta(s)=\sum_{n=1}^\infty \frac{1}{n^s},\;{\rm Re}(s)>1.
\end{equation*}
The series is analytic everywhere except for a simple pole at $s=1$ with residue 1, and the function, for ${\rm Re}(s) > 1$
\begin{equation*}
    \xi(s)=\frac{1}{2}s(s-1)\pi^{-s/2}\Gamma\left(\frac{s}{2}\right)\zeta(s)
\end{equation*}
is the restriction to $({\rm Re}(s) > 1)$ of a unique entire analytic function $\xi$, with satisfies the functional equation.
\begin{equation*}
    \xi(s)=\xi(1-s),\;s \in \mathbb{C}.
\end{equation*}
Our main purpose in this paper is to discuss a new class of functional equations by using both Poisson’s summation formula and Jacobi’s type theta function. The class of Riemann’s type functional equations are derived from self-reciprocal probability density functions, which is given in section 2, the second order of Zeta functional
equations for Riemann’s type is also investigated.
\section{Main Results}
In this article, we use the following notations, the probability density functions are denoted by $f$, $g$, $\cdots$ etc. and their characteristic functions $\hat{f}$, $\hat{g}$, $\cdots$ etc,. respectively.\\\\
There is an important formula, called Poisson’s summation formula, which has many
applications. The formula can be expressed in different ways. The following form of the formula is very convenient for us.
\begin{Lemma}[Poisson summation formula]\label{lemma1} Let $\hat{f}$ be a real symmetric characteristic function such that $\int_{\mathbb{R}}{|\hat{f}(t)|}dt<\infty$ and assume that $\hat{f}$ is non-increasing on $[0, \infty)$.
Then we have
$$\sum_{n=-\infty}^\infty \hat{f}(n)=\sum_{n=-\infty}^\infty \int_{\mathbb{R}}{\hat{f}(t)e^{-2\pi n i t}}dt,$$each series being absolutely convergent.
\end{Lemma}
\begin{Proof}
To see: T.M. Apostol, Mathematical Analysis, \nth{2}
Edition, P.332, Theorem 11.24. (1974).\hfill\qed
\end{Proof}\\\\
Note that a symmetric probability density function (pdf) $f$ on $\mathbb{R}$ is said to be selfreciprocal (SR) if
\begin{equation*}
    \hat{f}(t)=\sqrt{2\pi}f(t),\;t \in \mathbb{R}
\end{equation*}
where $\hat{f}$ is the characteristic function of $f$, that is.
\begin{equation*}
    \hat{f}(t)=\int_{\mathbb{R}}{e^{itx}}f(x)dx,\;t \in \mathbb{R}.
\end{equation*}
This curious property (that its density function and characteristic function are
essentially the same functions) is satisfied by two important functions, that are $N(0,1)$-distribution and hyperbolic cosh-distribution.\\
Fox examples, the $N(0,1)$-denisty $f_N(x)=\frac{1}{\sqrt{2\pi}} {e^{-x^2/2}},\;x \in \mathbb{R},$ is self-reciprocal, its ch.$f$. is given by $\hat{f}_N (x)=e^{-x^2/2},\;x \in \mathbb{R};$ the hyperbolic cosh-distribution with the pdf
$$f_c(x)=\frac{1}{\sqrt{2\pi}\cosh(\sqrt{\frac{\pi}{2}}x)},\;x \in \mathbb{R}$$
is also SR, its ch.$f$. is given by
$$\hat{f}_c(x)=\frac{1}{\cosh(\sqrt{\frac{\pi}{2}}x)},\;x \in \mathbb{R}$$
Note that the class of such self-reciprocal densities is surprisingly large (To see, L, Bondesson (1992), P.122).\\\\
The following lemma gives the most important sub class of self-reciprocal densities,
which is related to a variance mixture of the $N(0, \sigma^2)$-density. Similar results can be found in Teugels(1971).\\
\begin{Lemma}[Bondesson, 1992, P.122]\label{lemma2}
Let $f$ be a variance mixture of $N(0, \sigma^2)$-density
\begin{equation*}
    f(x)=\int_{0}^{\infty} \frac{1}{\sqrt{2\pi y}} e^{-\frac{x^2}{2y}}g(y)dy,\;x \in  \mathbb{R}
\end{equation*}
wheree $g$ is the mixing density for $\sigma^2 > 0$; assume that $g$ satisfies the relationship
\begin{equation*}
    g(x)=x^{-\frac{3}{2}}g\Big(\frac{1}{x}\Big),\;x >0.
\end{equation*}
Then $f$ is a self-reciprocal density, that is,
\begin{equation*}
    \hat{f}(x)=\sqrt{2\pi}f(x),\;x \in  \mathbb{R}.
\end{equation*}
\end{Lemma}
\begin{Proof}
This lemma follows from the unigueness of the Laplace transform (LT), Note that
the characteristic function of $N(0, \sigma^2)$-density is given by
    $$e^{-\frac{1}{2}\sigma^2t^2}=\int_{\mathbb{R}}e^{itx}\frac{1}{\sqrt{2\pi}\sigma}e^{-\frac{x^2}{2\sigma^2}}dx,\;t \in \mathbb{R}.$$\hfill\qed
\end{Proof}
\begin{flushleft}Note that the Lemma \ref{lemma2} can be rewrited by the Brownian subordinator (To see, Pitman and Yor, 2003) as follows.
Let $\beta_t,\;t\geq 0,$ be a standard Brownian motion, that is, the L\'evy process such that $\beta_t$ has normal distribution with $E(\beta_t)=0$ and $E(\beta_{t}^{2})=t.$ Let random variable $Y$ obey the mixing density $g$ as given in Lemma \ref{lemma2}, and $g$ satisfy the relationship
\begin{equation*}
    g(x)=x^{-\frac{3}{2}}g\left(\frac{1}{x}\right),\;x>0.
\end{equation*}
If $\beta_t,\;t\geq 0,$ is independent of the subordinator $Y$ then the random variable $\beta_Y$ has the self-reciprocal property, that is,
\begin{equation*}
    \hat{f}_{\beta_{Y}}(x)=\sqrt{2\pi}f_{\beta_{Y}}(x),\;x \in \mathbb{R}
\end{equation*}
\end{flushleft}
\definition[Poisson's summation property]\label{def1}
let $f$ be a symmetric probability density function $(p.d.f)$ on $\mathbb{R}$, and $\hat{f}$ be its characteristic function. The pdf $f$ $(or\;ch.f.\hat{f})$ is said to have the Poisson's summation property, if $f$ satisfies the following two conditions:
\begin{enumerate}
\item[$(1)$]
The symmetric ch.$f$.$\hat{f}$ is decreasing on $[0, \infty)$ and such that 
$\int_{\mathbb{R}}|\hat{f}(t)|dt<\infty.$
\item[$(2)$]
The pdf $f$ is self-reciprocal (SR), that is, $f$ satisfies $\hat{f}(t)=\sqrt{2\pi}f(t),\;t \in \mathbb{R}.$
\end{enumerate}
\begin{Lemma}\label{lemma3}
    Let $f$ be a symmetric pdf on $\mathbb{R}$, and $\hat{f}$ be the corresponding ch.$f.$, that is,
    $$\hat{f}(t)=\int_{\mathbb{R}} e^{itx}f(x)dx,\;t \in \mathbb{R}.$$
    Assume that $f$ has the Poisson's summation property as given in Definition \ref{def1}. Then, the following two identities hold,
    $$\sum_{n=-\infty}^{\infty}f(\sqrt{2\pi}xn)=\frac{1}{x}\sum_{n=-\infty}^{\infty}f\left(\frac{\sqrt{2\pi}n}{x}\right),\;x>0$$
    and
    $$\sum_{n=-\infty}^{\infty}\hat{f}(\sqrt{2\pi}xn)=\frac{1}{x}\sum_{n=-\infty}^{\infty}\hat{f}\left(\frac{\sqrt{2\pi}n}{x}\right),\;x>0.$$
\end{Lemma}

\begin{Proof}
\rm{The pdf $f$ is symmetric, and the first condition of Definition \ref{def1} implies that the pdf $f$ (continuous and bounded) is given by (the inversion theorem)
$$f(x)=\frac{1}{2\pi}\int_{\mathbb{R}}e^{-itx}\hat{f}(t)dt,\;x \in \mathbb{R}$$
and the condition of lemma \ref{lemma1} is also satisfied; hence we get, by using the Poisson summation formula, since $f$ is symmetric,
    \begin{eqnarray*}
    \sum_{n=-\infty}^{\infty}\hat{f}(n)&=&\sum_{n=-\infty}^{\infty}\int_{\mathbb{R}}\hat{f}(t)e^{-2\pi nit}dt\\&=&\sum_{n=-\infty}^{\infty}2\pi f(2\pi n).
    \end{eqnarray*} 
    Now by a suitable change of scale, and the second condition of Definition \ref{def1} implies the two identities hold. The proof is complete.\hfill\qed}
\end{Proof}
\begin{definition}
[Jacobi's Type theta function]\label{def2}
let $f$ be a symmetric probability density function and $\hat{f}$ its characteristic function as given in Definition \ref{def1}. The Jacobi's type theta function is defined for all $x>0$ by the equation
\begin{equation*}
    \theta(x)=\sum_{n=-\infty}^{\infty} \hat{f}(\sqrt{2 \pi}xn).
\end{equation*}
It follows from Lemma \ref{lemma3} that the theta function satisfies
\begin{equation*}
    \theta(x)=\frac{1}{x}\theta\big(\frac{1}{x}\big),\;x>0
\end{equation*}
and
\begin{equation*}
    \theta(x)=2\psi(x)+1,\;x>0
\end{equation*}
where the function $\psi$ is defined by
\begin{equation*}
\psi(x)=\sum_{n=1}^{\infty}\hat{f}(\sqrt{2\pi}xn),\;x>0.
\end{equation*}
\end{definition}
    \rm{The following Theorem \ref{thm1} gives a new class of Riemann's type functional equations by using both the Poisson's summation property and the Mellin's transform.}
\begin{theorem}[Riemann's Type functional equation]\label{thm1}
Let $f$ be a symmetric pdf on ${\mathbb{R}}$ and have the Poisson's summation property as given in Definition \ref{def1}. Assume that the moments of $f$ are finite for all orders. Then the Mellin's transform $\eta(s)$ of the $\psi$ function is equal to, for $s>1$,
\begin{eqnarray*}
    \eta(s)&\equiv&\int_{0}^{\infty}x^{s-1}\psi(x)dx\\
    &=&\left({\frac{1}{\sqrt{2\pi}}}\right)^s\zeta(s)\int_{0}^{\infty}x^{s-1}\hat{f}(x)dx
\end{eqnarray*}
    or
\begin{eqnarray*}
    \eta(s)=\left({\frac{1}{\sqrt{2\pi}}}\right)^{s-1}\zeta(s)\int_{0}^{\infty}x^{s-1}{f}(x)dx,\;s>1
\end{eqnarray*}
and the function $\xi(s)\equiv s(s-1)\eta(s)$ is an entire function of $s \in \mathbb{C}$ and satisfies the functional equation$\xi(s)=\xi(1-s)$.
\end{theorem}
\begin{Proof}
\rm{The symmetric pdf $f$ on $\mathbb{R}$ has the Poisson's summation property, it follows that the Mellin's transform $\eta(s)$ of the $\psi$ function is, for $s>1,$}
\begin{eqnarray*}
    \eta(s)&\equiv& \int_{0}^{\infty}x^{s-1}\psi(x)dx\\
    &=&\int_{0}^{\infty}x^{s-1}\sum_{n=1}^{\infty}\hat{f}(\sqrt{2\pi}xn)dx\\
    &=&\sum_{n=1}^{\infty}\int_{0}^{\infty}x^{s-1}\hat{f}(\sqrt{2\pi}xn)dx\\
    &=&\sum_{n=1}^{\infty}\int_{0}^{\infty}\big(\frac{y}{\sqrt{2\pi}n}\big)^{s-1}\hat{f}(y)\frac{1}{\sqrt{2\pi}n}dy\\
    &=&\sum_{n=1}^{\infty}\big(\frac{1}{\sqrt{2\pi}n}\big)^s\int_{0}^{\infty}y^{s-1}\hat{f}(y)dy\\
    &=&\big(\frac{1}{\sqrt{2\pi}}\big)^s\zeta(s)\int_{0}^{\infty}x^{s-1}\hat{f}(x)dx\\
    &=&\big(\frac{1}{\sqrt{2\pi}}\big)^{s-1}\zeta(s)\int_{0}^{\infty}x^{s-1}{f}(x)dx<\infty\\
\end{eqnarray*}
the last equality follows from the self-reciprocal probability and the moments of $f$ are finite for all orders.\\
\rm{On the other hand by the definition \ref{def2} of the Jacobi's theta function, it follows that the $\psi$ function satisfies}
$$\psi(x)=\frac{1}{x}\psi\left(\frac{1}{x}\right)+\frac{1}{2x}-\frac{1}{2},\;x>0.$$
and we have, for $s>1,$
\begin{eqnarray*}
    &&\int_{0}^{\infty}x^{s-1}\psi(x)dx\\
    &=&\int_{0}^{1}x^{s-1}\psi(x)dx+\int_{1}^{\infty}x^{s-1}\psi(x)dx\\
    &=&\int_{1}^{\infty}x^{s-1}\psi(x)dx+\int_{0}^{1}x^{s-1}\left[\frac{1}{x}\psi\left(\frac{1}{x}\right)+\frac{1}{2x}-\frac{1}{2}\right]dx\\
    &=&\int_{1}^{\infty}x^{s-1}\psi(x)dx+\int_{0}^{1}x^{s-2}\psi\left(\frac{1}{x}\right)dx+\int_{0}^{1}\frac{1}{2}x^{s-2}dx-\int_{0}^{1}\frac{1}{2}x^{s-1}dx\\
    &=&\int_{1}^{\infty}x^{s-1}\psi(x)dx+\int_{1}^{\infty}\left(\frac{1}{y}\right)^{s-2}\psi\left(y\right)\frac{1}{y^2}dy+\frac{1}{2(s-1)}-\frac{1}{2}\frac{1}{s}\\
    &=&\int_{1}^{\infty}x^{s-1}\psi(x)dx+\int_{1}^{\infty}x^{-s}\psi(x)dx+\frac{1}{2s(s-1)}\\
    &=&\int_{1}^{\infty}[x^{s-1}+x^{-s}]\psi(x)dx+\frac{1}{2s(s-1)}.
\end{eqnarray*}
Hence, we obtain that the $\eta$ function satisfies the functional equation $\eta(s)=\eta(1-s).$\\\\
Note that the symmetric pdf $f$ on ${\mathbb{R}}$ having moments of all orders, it implies that the Mellin's transform $\int_{0}^{\infty}x^{s-1}f(x)dx$ converges absolutely for $\rm{Re}(s)\geq 1.$ The convergence is also uniform in every half-plane $\rm{Re}(s)\geq 1+\epsilon,\;\epsilon>0,$ so this Mellin's transform is an analytic function of $s$ in the half-plane $\rm{Re}(s)\geq 1.$\\
Finally, we note that the $\eta(s)$ function has simple poles at$s=0$ and $s=1.$ Following Riemann, we multiply $\eta(s)$ by $s(s-1)$ to remove the poles and defined the Xi function $\xi(s)=s(s-1)\eta(s)$ and the result holds for all $s \in \mathbb{C}$ by analytic continuation. The proof is complete.\hfill\qed
\end{Proof}\\
\\
\rm{The pdf $f$ is symmetric, it follows that $\xi(0)=\frac{1}{2}=\xi(1)$.} Since Theorem \ref{thm1} implies
\begin{eqnarray*}
 &&\lim_{s \rightarrow 1 }\xi(s)\\
 &=&\lim_{s \rightarrow 1 }s(s-1)\eta(s)\\
 &=&\lim_{s \rightarrow 1 }s(s-1)\left(\frac{1}{\sqrt{2\pi}}\right)^{s-1}\zeta(s)\int_{0}^{\infty}x^{s-1}f(x)dx\\
  &=&\int_{0}^{\infty}f(x)dx\\
  &=&\frac{1}{2}
\end{eqnarray*}
and the result follows.\\
The first application of Theorem \ref{thm1} is the functional equation for the Riemann Zeta function as follows.
\begin{corollaryenv}$($proofs in Appendix$)$
Let $f_{N}$ be the  $N\left(0,1\right)$-density,
\begin{equation*}
    f_{N}(x)=\frac{1}{\sqrt{2\pi}}e^{-x^2/2},\;x \in \mathbb{R}.
\end{equation*}
Then the Mellin's transform $\eta_{N}(s)$ of the $\psi_{N}$ function, where
\begin{equation*}
    \psi_{N}(x)=\sum_{n=1}^{\infty} \hat{f}_{N}(\sqrt{2\pi}xn),\;x>0,
\end{equation*}
is
\begin{eqnarray*}
    \eta_{N}(s)&\equiv&\int_{0}^{\infty}x^{s-1}\psi_{N}(x)dx\\
               &=&\frac{1}{2}\pi^{-s/2}\Gamma\left(\frac{s}{2}\right)\zeta(s),\;s>1
 \end{eqnarray*}
and the Xi function, for $s \in \mathbb{C}$, is
\begin{eqnarray*}
    \xi_{N}(s)&\equiv&s(s-1)\eta_{N}(s)\\
              &=&\frac{1}{2}s(s-1)\pi^{-s/2}\Gamma\left(\frac{s}{2}\right)\zeta(s),
\end{eqnarray*}
the $\xi_{N}(s)$ is an entire function of $s$ and satisfies the functional equation
$\xi_{N}(s)=\xi_{N}(1-s).$
\end{corollaryenv}

\begin{corollaryenv}$($proofs in Appendix$)$
Let $f_{c}$ be the pdf of the hyperbolic cosh-distribution on $\mathbb{R}$, that is,
\begin{equation*}
    f_{c}(x)=\frac{1}{\sqrt{2\pi}\cosh{\left(\sqrt{\frac{\pi}{2}}x\right)}},\;x \in \mathbb{R}.
\end{equation*}
Then the Mellin's transform $\eta_{c}(s)$ of the $\psi_{c}$ function, where
\begin{equation*}
    \psi_{c}(x)=\sum_{n=1}^{\infty}\hat{f}_{c}(\sqrt{2\pi}xn),\;x>0,
\end{equation*}
is
\begin{eqnarray*}
    \eta_{c}(s)&\equiv&\int_{0}^{\infty}x^{s-1}\psi_{c}(x)dx\\
               &=&2\cdot\pi^{-s}\Gamma(s)\zeta(s)\sum_{n=0}^{\infty}\frac{(-1)^{n}}{(2n+1)^{s}},\;s>1,
\end{eqnarray*}
and the Xi function for $s \in \mathbb{C}$, is
\begin{eqnarray*}
    \xi_{c}(s)&\equiv&s(s-1)\eta_{c}(s)\\
              &=&2s(s-1)\pi^{-s}\Gamma(s)\zeta(s)\sum_{n=0}^{\infty}\frac{(-1)^{n}}{(2n+1)^{s}},
\end{eqnarray*}
the $\xi_{c}(s)$ is an entire function of $s$ and satisfies the functional equation $\xi_{c}(s)=\xi_{c}(1-s).$
\end{corollaryenv}
Note that Lemma \ref{lemma2} above gives the most important subclass of self-reciprocal density. Now by using Theorem \ref{thm1} and following Riemann's technique, we can obtain a class of the second order Zeta functional equations for Riemann's type, as follows.
\begin{theorem}[Second order case]\label{thm2}
Let pdf $f_M$ be a variance mixture of $N\left(0, \sigma^2\right)$-density,
$$f_M(x)=\int_0^{\infty} \frac{1}{\sqrt{2 \pi y}} e^{-\frac{x^2}{2 y}} \cdot g(y) d y, x \in \mathbb{R},$$
Where $g$ is the mixing density for $\sigma^2>0$, and a function $\xi_{1}$ be defined by, for $s>1$
$$
\xi_1(s)=\int_0^{\infty} x^{\frac{s-1}{2}} g(x) d x.
$$
If the mixing density g satisfies the relationship
$$
g(x)=x^{-\frac{3}{2}} \cdot g\left(\frac{1}{x}\right), \quad x>0
$$
and the moments of $g$ are finite for all orders, then the Mellin's transform $\eta_M(s)$ of the $\psi_M$ function,
$$
\psi_M(x)=\sum_{n=1}^{\infty} \hat{f}_M(\sqrt{2 \pi} x n), \quad x>0,
$$
is equal to, for $s>1$
$$
\begin{aligned}
\eta_M(s) & \equiv \int_0^{\infty} x^{s-1} \psi_M(x) d x \\
&=\frac{1}{2} \pi^{-s / 2} \Gamma\left(\frac{s}{2}\right) \zeta(s) \cdot \int_0^{\infty} x^{\frac{s-1}{2}} g(x) d x
\end{aligned}
$$
and the function, for $s>1$,
$$
\xi_2(s) \equiv s(s-1) \eta_M(s)
$$
is equal to a product of two entire functions, that is.
$$
\xi_2(s)=\xi_0(s) \xi_1(s), \quad s \in \mathbb{C}
$$
where $\xi_0(s) \equiv \xi(s)$ is the functional equation for the Riemann Zeta function,
$$
\xi_0(s)=\frac{1}{2} s(s-1) \pi^{-s / 2} \Gamma(\frac{s}{2}) \zeta(s)
$$
and
$$
\xi_1(s) \equiv \int_0^{\infty} x^{\frac{s-1}{2}} g(x) d x, \quad s>1.
$$
Furthermore, all the three functions $\xi_1(s)$, $i=0,1,2$, are entire functions of $s$ and satisfy the functional equation
$$
\xi_i(s)=\xi_i(1-s),
$$
\end{theorem}

\begin{Proof}
Firstly, note that the mixing density $g$ satisfies the relationship
$$
g(x)=x^{\frac{-3}{2}} g\left(\frac{1}{x}\right), x>0,
$$
which is equivalent to
$$
x^{\frac{3}{4}} g(x)=x^{-\frac{3}{4}} g\left(\frac{1}{x}\right), x>0
$$
by using Lemma 2 , it follows that the pdf $f_M$ is a self-reciprocal density, that is,
$$
\hat{f}_M(x)=\sqrt{2 \pi} f_M(x), \quad x \in \mathbb{R} .
$$
by Definition \ref{def1}, it follows that $f_M$ has the Poisson's summation property. Since the moments of $g$ are finite at for all orders, and so are $f_M$. Now by using Theorem \ref{thm1}, it follows that the Mellin's transform $\eta_M(s)$ of the $\psi_M$ function is equal to, for $s>1$,
$$
\begin{aligned} \eta_M(s)\equiv & \int_0^{\infty} x^{s-1} \psi_M(x) d x \\ = & \left(\frac{1}{\sqrt{2 \pi}}\right)^{s-1} \zeta(s) \cdot \int_0^{\infty} x^{s-1} f_M(x) d x \\ = & \left(\frac{1}{\sqrt{2 \pi}}\right)^{s-1} \zeta(s) \int_0^{\infty} x^{s-1} \int_0^{\infty} \frac{1}{\sqrt{2 \pi y}} e^{-\frac{x^2}{2 y}} g(y) d y d x \\ = & \left(\frac{1}{\sqrt{2 \pi}}\right)^{s-1} \zeta(s) \int_0^{\infty} \frac{1}{\sqrt{2 \pi}} y^{-1 / 2} g(y) \cdot \int_0^{\infty} x^{s-1} e^{-\frac{x^2}{2 y}} d x d y \\ = & \left(\frac{1}{\sqrt{2 \pi}}\right)^s \zeta(s) \int_0^{\infty} y^{\frac{-1}{2}} \cdot g(y) \cdot \frac{\Gamma(\frac{s}{2})}{2\left(\frac{1}{2 y}\right)^{s / 2}} d y \\ = & \frac{1}{2} \pi^{-s / 2} \Gamma\left(\frac{s}{2}\right) \zeta(s) \cdot \int_0^{\infty} y^{\frac{s-1}{2}} g(y) d y
\end{aligned}$$
%and the $X_i$ function is
%$$
%\xi_M(s) \equiv s(s-1) \eta_M(s),
%$$
%The $\xi_M$ is an entire function of $s \in \mathbb{C}$ and satisfies the functional equation
%$$
%\xi_M(s)=\xi_M(1-s).
%$$
%The proof is complete.\qed\\
and the function $\xi_2$ is equal to
$$
\begin{aligned}
\xi_2(s) & \equiv s(s-1) \eta_M(s) \\
& =\frac{1}{2} s(s-1) \pi^{-s / 2} \Gamma\left(\frac{s}{2}\right) \zeta(s) \int_0^{\infty} x^{\frac{s-1}{2}} g(x) d x \\
& \equiv \xi_0(s) \cdot \xi_1(s)
\end{aligned}
$$
$\xi_2$ is an entire function of $s$ and satisfies the functional equation
$$
\xi_z(s)=\xi_2(1-s) \text {. }
$$
For completeness, following Riemann, we can prove that the function $\xi_1$ also satisfies the functional equation
$$
\xi_1(s)=\xi_1(1-s) \text {. }
$$
Since the mixing density $g$ satisfies the relationship
$$
g(x)=x^{-\frac{3}{2}} g\left(\frac{1}{x}\right), \quad x>0
$$
it follows that
$$
\begin{aligned}
\xi_1(s) & =\int_0^{\infty} x^{\frac{s-1}{2}} g(x) d x \\
& =\int_0^1 x^{\frac{s-1}{2}} g(x) d x+\int_1^{\infty} x^{\frac{s-1}{2}} g(x) d x \quad \\
& =\int_1^{\infty}\left(\frac{1}{y}\right)^{\frac{s-1}{2}} \cdot g\left(\frac{1}{y}\right) \frac{1}{y^2} d y+\int_1^{\infty} x^{\frac{s-1}{2}} g(x) d x\;\;(\because\text{set}\;\;x=\frac{1}{y}) \\
& =\int_1^{\infty} y^{-\frac{s}{2}} \cdot y^{-\frac{3}{2}} \cdot g\left(\frac{1}{y}\right) d y+\int_1^{\infty} x^{\frac{s-1}{2}} g(x) d x \\
& =\int_1^{\infty} x^{-\frac{s}{2}} \cdot g(x) d x+\int_1^{\infty} x^{\frac{s-1}{2}} g(x) d x \\
& =\int_1^{\infty}\left(x^{-\frac{s}{2}}+x^{\frac{s-1}{2}}\right) g(x) d x \\
& =\xi_1(1-s)
\end{aligned}
$$
the other results follow immediately and the proof is complete.\hfill\qed
\end{Proof}
\\
Note that the Mellin's transform $\eta_M(s)$ as given in Theorem \ref{thm2} has simple poles at $s=0$ and $s=1$. Following Riemann, we multiply by $s(s-1)$ to remove the poles.
It follows that
$\xi_2(0)=\frac{1}{2}=\xi_2(1)$
and $\xi_1(0)=1=\xi_1(1)$. Since
$$
\begin{aligned}
\lim _{s \rightarrow 1} \xi_2(s) & =\lim _{s \rightarrow 1} s(s-1) \eta_M(s) \\
& =\lim _{s \rightarrow 1} \xi(s) \cdot \int_0^{\infty} x^{\frac{s-1}{2}} g(x) d x \\
& =\frac{1}{2} \quad\left(\because g \text { is a pdf}  \text { on } \mathbb{R}_{+}\right) .
\end{aligned}
$$
%%%%%%%%%%%%%%%%%%%%%%%%%%%%%%%%%%%%%%%%%%%%%%%%%%%%%%%%%%%%%%%%%%%%%
%%%%%%%%%%%%%%%%%%%%%%%%%%%%%%%%%%%%%%%%%%%%%%%%%%%%%%%%%%%%%%%%%%%%%
%%%%%%%%%%%%%%%%%%%%%%%%%%%%%%%%%%%%%%%%%%%%%%%%%%%%%%%%%%%%%%%%%%%%%
\begin{corollaryenv}$($proofs in Appendix$)${$($Sinh-density $f_{S_2}$ case$)$}

\noindent Let random Variable $\hat{z}$ obey a distribution with pdf $f_{\hat{z}}$ and $f_{\hat{z}}$ be a variance mixture of $N\left(0,\sigma^2\right)$-density,
$$
f_{\hat{z}}(x)=\int_0^{\infty} \frac{1}{\sqrt{2 \pi y}} e^{-\frac{x^2}{2 y}} \cdot f_z(y) d y, x \in \mathbb{R}
$$
where $z \geqslant 0$ is a non-negative random variabl with the mixing density $f_z$,
$$
f_z(x)=2 \sqrt{x} \sum_{n=1}^{\infty}\left(2 \pi^2 n^4 x-3 \pi n^2\right) e^{-\pi n^2 x}, x>0
$$
Then, the following three statements hold.
\begin{description}
\item[(a).] If function $\xi_z(s)$ is defined by, for $s>1$,
$$
\xi_{z}(s)=\int_0^{\infty} x^{\frac{s-1}{2}} f_{z}(x) d x
$$
then
$$
\xi_z(s)=2 \xi(s), \quad s \in \mathbb{C}
$$
Where $\xi$ is the functional equation for the Riemann Zeta function, that is,
$$
\xi(s)=\frac{1}{2} s(s-1) \pi^{-s / 2} \cdot \Gamma\left(\frac{s}{2}\right) \cdot \zeta(s) .
$$
\item[(b).] If $\eta_{\hat{z}}(s)$ is the Mellin's transform of the $\psi_{\hat{z}}$ function,
$$
\psi_{\hat{z}}(x)=\sum_{n=1}^{\infty} \hat{f}_{\hat{z}}(\sqrt{2 \pi} x n), \quad x>0
$$
then
$$
\begin{aligned}
\eta_{\hat{z}}(s) & \equiv \int_0^{\infty} x^{s-1} \psi_{\hat{z}}(x) d x \\
& =\frac{1}{2} \pi^{-s / 2} \Gamma\left(\frac{s}{2}\right) \zeta(s) \cdot \xi_z(s) .
\end{aligned}
$$
\item[(c).]  Furthermore, the function
$$
\xi_{\hat{z}}(s) \equiv s(s-1) \eta_{\hat{z}}(s)
$$
is equal to
$$
\xi_{\hat{z}}(s)=2 \cdot[\xi(s)]^2, \quad s \in \mathbb{C}
$$
\end{description}
\end{corollaryenv}
\begin{corollaryenv}$($proofs in Appendix$)$
Let $f_{\hat{w}}$ be a variance mixture of $N\left(0,\sigma^2\right)^2$-density
$$
f_{\hat{w}}(x)=\int_0^{\infty} \frac{1}{\sqrt{2 \pi y}} e^{-\frac{x^2}{2 y}} f_w(y) d y,\quad x \in \mathbb{R}
$$
where $f_w$ is the mixing density
$$
f_w(x)=\frac{1}{2 \sqrt{x}} H_1(\sqrt{x}), \quad x>0
$$
and the function $H_1(x)$ is, for $x>0$
$$
H_1(x)=2 \cdot \sum_{x=1}^{\infty} n \pi x \operatorname{sech}(n \pi x)\left[n \pi x\left(\tanh ^2(n \pi x)-\operatorname{sech}^2(n \pi x)\right)-2 \tanh (n \pi x)\right] \text {. }
$$
Then the following three statements hold.
\begin{description}
\item[(a).] If function $\xi_w$ is defined by, for $s>1$
$$
\xi_w(s)=\int_0^{\infty} x^{\frac{s-1}{2}} f_w(x) d x
$$
then
$$
\begin{aligned}
\xi_w(s) & \equiv 2 \xi_c(s) \text { (as given in corollary 2) } \\
& =4 s(s-1) \pi^{-s} \Gamma(s) \zeta(s) \sum_{n=0}^{\infty} \frac{(-1)^n}{(2 n+1)^s}
\end{aligned}
$$
$\xi_w$ is an entire function of $s \in \mathbb{C}$ and satisfies the functional equation
$$
\xi w(s)=\xi_w(1-s) .
$$
\item[(b).] If $\eta_{\hat{w}}(s)$ is the Mellin's transform of the $\psi_{\hat{w}}$ function.
$$
\psi_{\hat{w}}(x)=\sum_{n=1}^{\infty} \hat{f}_{\hat{w}}(\sqrt{2 \pi} x n), \quad x>0
$$
then
$$
\begin{aligned}
\eta_{\hat{w}}(s) & \equiv \int_0^{\infty} x^{s-1} \psi_{\hat{w}}(x) d x \\
& =\frac{1}{2} \pi^{-s / 2} \Gamma\left(\frac{s}{2}\right) \zeta(s) \cdot \xi_w(s)
\end{aligned}
$$
\item[(c).] Furthermore, the function
$$
\xi_{\hat{w}}(s) \equiv s(s-1) \eta_{\hat{w}}(s)
$$
is equal to
$$
\xi_{\hat{w}}(s)=\xi(s) \cdot \xi_w(s)
$$
$\xi_{\hat{w}}$ is an entire function of $s \in \mathbb{C}$ and satisfies the functional equation
$$
\xi_{\hat{w}}(s)=\xi_{\hat{w}}(1-s)
$$
\end{description}
\end{corollaryenv}
\begin{corollaryenv}$($proofs in Appendix$)${(cosh-density $f_{C_1}$ case)}
Let random variable $\hat{T}$ obey a distribution with pdf $f_{\hat{T}}$ and $f_{\hat{T}}$ be a variance mixture of the $N\left(0, \sigma^2\right)$-density,
$$
f_{\hat{T}}(x)=\int_0^{\infty} \frac{1}{\sqrt{2 \pi y}} e^{-\frac{x^2}{2 y}} \cdot f_T(y) d y, \quad x \in \mathbb{R}
$$
where $T \geqslant 0$ is a non-negative random variable with the mixing density $f_T$
$$
f_T(x)=2 \sum_{n=0}^{\infty}(-1)^n\left(n+\frac{1}{2}\right) \cdot e^{-\left(n+\frac{1}{2}\right)^2 \pi x},\quad x>0 .
$$
Then the following three statements hold.
\begin{description}
\item[(a).] If function $\xi_4(s)$ is defined by, for $s>1$,
$$
\xi_4(s)=\int_0^{\infty} x^{\frac{s-1}{2}} f_T(x) d x
$$
then 
$$
\xi_4(s)=\left(\frac{4}{\pi}\right)^{\frac{s+1}{2}} \Gamma\left(\frac{s+1}{2}\right) \cdot L_{x_4}(s)
$$
where
$$
L_{x_4}(s)=\sum_{n=0}^{\infty} \frac{(-1)^n}{(2 n+1)^s}, \quad R_e(s)>0
$$
The function $\xi_4$ satisfies the functional equation
$$
\xi_4(s)=\xi_4(1-s) \text {, }
$$
and $\xi_4$ is an entire function of $s \in \mathbb{C}$.
\item[(b).] If $\eta_{\hat{T}}{ }(s)$ is the Mellin's transform of the $\psi_{\hat{T}}$ function,
$$
\psi_{\hat{T}}(x)=\sum_{n=1}^{\infty} \hat{f}_{\hat{T}}(\sqrt{2 \pi}xn), \quad x>0
$$
then
$$
\begin{aligned}
\eta_{\hat{T}}(s) & \equiv \int_0^{\infty} x^{s-1} \psi_{\hat{T}}(x) d x \\
& =\frac{1}{2} \pi^{-s / 2} \Gamma\left(\frac{s}{2}\right) \zeta(s) \xi_4(s)
\end{aligned}
$$
\item[(c).] Furthermore, the function
$$
\xi_{\hat{T}}(s) \equiv s(s-1) \eta_{\hat{T}}(s)
$$
is equal to
$$
\xi_{\hat{T}}(s)=\xi(s) \cdot \xi_4(s), \quad s \in \mathbb{C} .
$$
and $\xi_{\hat{T}}$ is an entire function of $s$.
\end{description}
\end{corollaryenv}
The class of GGC's was introduced by O. Thorin in 1977 as follows. (To see, Bondesson 1992)\\The following Theorem \ref{thm3} is concerned with the class GGC.
\begin{theorem}[Generalized Gamma Convolution Case]\label{thm3}
Let $f_1$ be a variance mixture of $N\left(0, \sigma^2\right)$-density 
$$
f_1(x)=\int_0^{\infty} \frac{1}{\sqrt{2 \pi y}} e^{-\frac{x^2}{2 y}} g_1(y) d y, \quad x \in \mathbb{R}
$$
where $g_1$ is the mixing density, if $g_1$ is of the following form
$$
g_1(x)=c \cdot x^{-3 / 4} h(x) \cdot h\left(\frac{1}{x}\right),\quad x>0
$$
where $c$ is a constant and $h$ is the
Laplace - Transform (LT) of a GGC,
and the moments of $g_1$ are finite for all
orders. Then the following three statements hold.
\begin{description}
\item[(a).] If function $\xi_{G_1}$ is defined by, for $s>1$
$$
\xi_{G_1}(s)=\int_0^{\infty} x^{\frac{s-1}{2}} g_1(x) d x
$$
then $\dis\xi_{G_1}$ satisfies the functional equation
$$
\xi_{G_1}(s)=\xi_{G_1}(1-s)
$$
and $\xi_{G_1}$ is an entire function of $s \in \mathbb{C}$.
\item[(b).] If $\eta_{\hat{G}_1}(s)$ is the Mellin's transform of the $\psi_{\hat{G}_1}$ function,
$$
\psi_{\hat{G}_1}(x)=\sum_{n=1}^{\infty} \hat{f}_1(\sqrt{2 \pi} x n),\quad x>0,
$$
then
$$
\begin{aligned}
\eta_{\hat{G}_1}(s) & \equiv \int_0^{\infty} x^{s-1} \psi_{\hat{G}_1}(x) d x \\
& =\frac{1}{2} \pi^{-\frac{s}{2}} \Gamma\left(\frac{s}{2}\right) \zeta(s) \xi_{G_1}(s).
\end{aligned}
$$
\item[(c).] Furthermore, the function
$$
\xi_{\hat{G}_1}(s) \equiv s(s-1) \eta_{\widehat{G}_1}(s)
$$
is equal to
$$
\xi_{\hat{G}_1}(s)=\xi(s) \cdot \xi_{G_1}(s), \quad s \in \mathbb{C}
$$
and $\xi_{\hat{G}_1}$ is an entire function of $s$, $\xi_{\hat{G}_1}$ satisfies the functional equation
$$
\xi_{\hat{G}_1}(s)=\xi_{\hat{G}_1}(1-s).
$$
\end{description}
\end{theorem}
\begin{Proof}
This Theorem follows from both Theorem 2 and Bondesson (1992, P.122 below). Since, it easily find that the pdf $f_1$ of a variance mixture with the mixing density $g_1 \in B$ is self-reciprocal iff
$$
g_1(x)=c \cdot x^{-\frac{3}{4}} \cdot h(x) \cdot h\left(\frac{1}{x}\right), \quad x>0
$$
Where $h$ is the Lapluce-transform of $a$ GGC. Note that
$$
\begin{aligned}
g_1(x) & =c \cdot x^{-\frac{3}{4}} h(x)h\left(\frac{1}{x}\right) \\
& =x^{-\frac{3}{2}} c \cdot\left(\frac{1}{x}\right)^{-\frac{3}{4}} \cdot h\left(\frac{1}{x}\right) \cdot h(x) \\
& =x^{-\frac{3}{2}} \cdot g_1\left(\frac{1}{x}\right), \quad x>0 .
\end{aligned}
$$
All the conditions of Theorem 2 hold and the results follow immediately. The proof is complete.\qed
\end{Proof}
\begin{corollaryenv}$($proofs in Appendix$)$
Let $f_2$ be a variance mixture of $N\left(0, \sigma^2\right)$-density
$$
f_2(x)=\int_0^{\infty} \frac{1}{\sqrt{2 \pi y}} e^{-\frac{x^2}{2 y}} \cdot g_2(y) d y,\quad x \in \mathbb{R}
$$
where $g_2$ is the mixing density, and
$$
g_2(x)=c \cdot x^{-\frac{3}{4}} e^{-a\left(x^\alpha+x^{-\alpha}\right)},\quad x>0
$$
where $0<\alpha \leqslant 1$ and $a>0, c>0$. Then the following three statements hold.
\begin{description}
\item[(a).]If function $\xi_{G_2}$ is defind by, for $s>1$
$$
\xi_{G_2}(s)=\int_0^{\infty} x^{\frac{s-1}{2}} \cdot g_2(x) d x
$$
then $\xi_{G_2}$ satisfies the functional equation
$$
\xi_{G_2}(s)=\xi_{G_2}(1-s)
$$
and $\xi_{G_2}$ is an entire function of $s \in \mathbb{C}$.
\item[(b).] If $\eta_{\hat{G}_2}(s)$ is the Mellin's transform of the $\psi_{\hat{G}_2}$ function,
$$
\psi_{\hat{G}_2}(x)=\sum_{n=1}^{\infty} \hat{f}_2(\sqrt{2 \pi} x n), \quad x>0
$$
then
$$
\begin{aligned}
\eta_{\hat{G}_2}(s) & \equiv \int_0^{\infty} x^{s-1} \psi_{\hat{G}_2}(x) d x \\
& =\frac{1}{2} \pi^{-\frac{s}{2}} \cdot \Gamma\left(\frac{s}{2}\right) \zeta(s) \cdot \xi_{G_2}(s).
\end{aligned}
$$
\item[(c).] Furthermore, the function
$$
\xi_{\hat{G}_2}(s) \equiv s(s-1) \eta_{\hat{G}_2}(s)
$$
is equal to
$$
\xi_{\hat{G}_2}(s)=\xi(s) \cdot \xi_{G_2}(s),\quad s \in \mathbb{C}
$$
and $\xi_{\hat{G}_2}$ is an entire function of $s$,
$\xi_{\hat{G}_2}$ satisfies the functional equation
$$
\xi_{\hat{G}_2}(s)=\xi_{\hat{G}_2}(1-s), \quad s \in \mathbb{C} .
$$
\end{description}
\end{corollaryenv}
\begin{corollaryenv}$($proofs in Appendix$)${$($Generlized Inverse Gaussian dist.$)$}
Let $f_3$ be a variance mixture of $N\left(0, \sigma^2\right)$-density
$$
f_3(x)=\int_0^{\infty} \frac{1}{\sqrt{2 \pi y}} e^{-\frac{x^2}{2 y}} \cdot g_3(y) d y, \quad x \in \mathbb{R}
$$
where $g_3$ is the mixing density, and
$$
g_3(x)=\frac{1}{2 K_{\frac{1}{4}
}{ }^{(a)}} x^{\frac{-3}{4}} e^{\frac{-a(x+\frac{1}{x})}{2}},\quad x>0
$$
where $K_p(a)$ is a modified Bessel function of the second kind and $a>0$ is a parameter. Then the following three statements hold.
\begin{description}
\item[(a).]If function $\xi_{G_3}$ is defined by, for $s>1$
$$
\xi_{G_3}(s)=\int_0^{\infty} x^{\frac{s-1}{2}} \cdot g_3(x) d x
$$
then 
$$
\xi_{G_3}(s)=\frac{K_{p_1}(a)}{K_{\frac{1}{4}}(a)}, \quad s \in \mathbb{C}
$$
where $p_1=\frac{2s-1}{4}$ and $a>0$, $\xi_{G_3}$ satisfies the functional equation
$$
\xi_{G_3}(s)=\xi_{G_3}(1-s), \quad s \in \mathbb{C}
$$
and $\xi_{G_3}$ is an entire function of $s$.
\item[(b).] If $\eta_{{\hat{G}_3}}$ (s) is the Mellin's transform of the $\psi_{\hat{G}_3}$ function,
$$
\psi_{\hat{G}_3}(x)=\sum_{n=1}^{\infty} \hat{f}_3(\sqrt{2 \pi} x n),\quad x>0
$$
then
$$
\begin{aligned}
\eta_{\hat{G}_3}(s) & \equiv \int_0^{\infty} x^{s-1} \psi_{\hat{G}_3}(x) d x \\
& =\frac{1}{2} \pi^{\frac{-s}{2}} \Gamma\left(\frac{s}{2}\right) \zeta(s) \cdot \xi_{G_3}(s)
\end{aligned}
$$
\item[(c).] Furthermore, the function
$$
\xi_{\hat{G}_3}(s) \equiv s(s-1) \eta_{\hat{G}_3}(s)
$$
is equal to
$$
\xi_{\hat{G}_3}(s)=\xi(s) \cdot \xi_{G_3}(s), \quad s \in \mathbb{C}
$$
and $\xi_{\hat{G}_3}$ is an entire function of $S$.
$\xi_{\hat{\epsilon}_3}$ satisfies the function equation
$$
\xi_{\hat{G}_3}(s)=\xi_{\hat{G}_3}(1-s), \quad s \in \mathbb{C}
$$
\end{description}
\end{corollaryenv}
\begin{corollaryenv}$($proofs in Appendix$)${$($Lévy formula case$)$}
Let $f_4$ be a variance mixture of $N\left(0, \sigma^2\right)$-density 
$$
f_4(x)=\int_0^{\infty} \frac{1}{\sqrt{2 \pi y}} e^{-\frac{x^2}{2 y}} g_4(y) d y,\quad x \in \mathbb{R}
$$
where $g_4$ is the mixing density, and $$ g_4(x)=\frac{1}{2}\left(\frac{\lambda}{\pi}\right)^{\frac{1}{2}} \cdot x^{\frac{-3}{4}} \cdot e^{-\lambda\left(x^{\frac{1}{4}}-x^{-\frac{1}{4}}\right)^2},\quad x>0 $$
where $\lambda>0$ is a parameter. The following three statements hold.
\begin{description}
\item[(a).] If function $\xi_{G_4}$ is defined by, for $s>1$
$$
\xi_{G_4}(s)=\int_0^{\infty} x^{\frac{s-1}{2}} \cdot g_4(x) d x \text {, }
$$
then $\xi_{G_4}$ satisfies the functional equation
$$
\xi_{G_4}(s)=\xi_{G_4}(1-s), \quad s \in \mathbb{C}
$$
and $\xi_{G_4}$ is an entire function of $s$.
\item[(b).] If $\eta_{\widehat{G}_4}(s)$ is the Mellin's transform of the $\psi_{\hat{G}_4}$ function.
$$
\psi_{\hat{G}_4}(x)=\sum_{n=1}^{\infty} \hat{f}_{4}(\sqrt{2 \pi}xn),\quad x>0
$$
then
$$
\begin{aligned}
\eta_{\hat{G_4}}(s) & \equiv \int_0^{\infty} x^{s-1} \psi_{\hat{G}_4}(x) d x \\
& =\frac{1}{2} \pi^{\frac{-s}{2}} \Gamma\left(\frac{s}{2}\right) \zeta(s) \cdot \xi_{G_4}(s).
\end{aligned}
$$
\item[(c).] Furthermore, the function
$$
\xi_{\hat{G}_4}(s) \equiv s(s-1) \eta_{\hat{G}_4}(s)
$$
is equal to
$$
\xi_{\hat{G}_4}(s)=\xi(s) \cdot \xi_{G_4}(s), \quad s \in \mathbb{C}
$$
and $\xi_{\hat{G}_4}$ is an entire function of $s$.
$\xi_{\hat{G}_4}$ satisfies the functional equation
$$
\xi_{\hat{G}_4}(s)=\xi_{\hat{G}_4}(1-s), \quad s \in \mathbb{C} .
$$
\end{description}
\end{corollaryenv}

\section{Appendix}
\begin{refproof}
    It is well-known that the $N(0, 1)-density$
    \begin{equation*}
        f_{N}(x)=\frac{1}{\sqrt{2\pi}}e^{-x^{2}/2},\;x \in \mathbb{R}
    \end{equation*}
satisfies all the conditions of Theorem \ref{thm1}. By using Theorem \ref{thm1}, it follows that the $\eta_{N}(s)$ function is given by, for $s>1$,
\begin{eqnarray*}
    \eta_{N}(s)&=&\left(\frac{1}{\sqrt{2\pi}}\right)^{s-1}\zeta(s)\int_{0}^{\infty}x^{s-1}f_{N}(x)dx\\
    &=&\left(\frac{1}{\sqrt{2\pi}}\right)^{s-1}\zeta(s)\int_{0}^{\infty}x^{s-1}\frac{1}{\sqrt{2\pi}}e^{-x^{2}/2}dx\\
    &=&\left(\frac{1}{\sqrt{2\pi}}\right)^{s}\zeta(s)\frac{\Gamma(s/2)}{2(1/2)^{s/2}}\\
    &=&\frac{1}{2}\pi^{-s/2}\Gamma\left(\frac{s}{2}\right)\zeta(s)
\end{eqnarray*}
and the Xi function, for $s \in \mathbb{C}$, is
\begin{eqnarray*}
    \xi_{N}(s)&\equiv&s(s-1)\eta_{N}(s)\\
    &=&\frac{1}{2}s(s-1)\pi^{-s/2}\Gamma\left(\frac{s}{2}\right)\zeta(s).
\end{eqnarray*}
The $\xi_{N}(s)$ is an entire function of $s$ and satisfies the functional equation 
$\xi_{N}(s)=\xi_{N}(1-s).$\\
$(To\;see,\;Apostol\;1976,\;p.260)$.\hfill\qed
\end{refproof}
\begin{refproof}
    It is also clearly that the hyperbolic cosh-distribution with the pdf
    \begin{equation*}
        f_{c}(x)=\frac{1}{\sqrt{2\pi}\cosh(\sqrt{\frac{\pi}{2}}x)},\;x \in \mathbb{R}
    \end{equation*}
    satisfies all the conditions of Theorem \ref{thm1}. By using Theorem \ref{thm1}, it follows that the $\eta_{c}$ function is given by, for $s>1$,
    \begin{eqnarray*}
    \eta_{c}(s)
    &=&\left(\frac{1}{\sqrt{2\pi}}\right)^{s-1}\zeta(s)\int_{0}^{\infty}x^{s-1}f_{c}(x)dx\\
    &=&\left(\frac{1}{\sqrt{2\pi}}\right)^{s-1}\zeta(s)\int_{0}^{\infty}x^{s-1}\cdot\frac{1}{\sqrt{2\pi}\cosh(\sqrt{\frac{\pi}{2}}x)}dx\\
    &=&\left(\frac{1}{{2\pi}}\right)^{s}\zeta(s)\int_{0}^{\infty}y^{s-1}\cdot\frac{1}{\cosh\left(\frac{y}{2}\right)}dy\\
    &=&\left(\frac{1}{{2\pi}}\right)^{s}\zeta(s)\int_{0}^{\infty}y^{s-1}\cdot\frac{2e^{-y/2}}{1+e^{-y}}dy\\
    &=&\left(\frac{1}{{2\pi}}\right)^{s}\zeta(s)\int_{0}^{\infty}y^{s-1}2\cdot e^{-y/2}\cdot\sum_{n=0}^{\infty}(-e^{-y})^{n}dy\\
    &=&2\left(\frac{1}{{2\pi}}\right)^{s}\zeta(s)\sum_{n=0}^{\infty}(-1)^{n}\int_{0}^{\infty}y^{s-1}\cdot e^{-(n+\frac{1}{2})y}dy\\
    &=&2\left(\frac{1}{{2\pi}}\right)^{s}\zeta(s)\sum_{n=0}^{\infty}(-1)^{n}\cdot\frac{\Gamma(s)}{(n+\frac{1}{2})^{s}}\\
    &=&2\cdot\pi^{-s}\Gamma(s)\zeta(s)\cdot\sum_{n=0}^{\infty}\frac{(-1)^{n}}{(2n+1)^{s}},\;s>1
\end{eqnarray*}
and the Xi function, for $s \in \mathbb{C}$.
\begin{eqnarray*}
    \xi_{c}(s)&\equiv&s(s-1)\eta_{c}(s)\\
    &=&2s(s-1)\pi^{-s}\Gamma(s)\zeta(s)\cdot \sum_{n=0}^{\infty}\frac{(-1)^{n}}{(2n+1)^{s}}.
\end{eqnarray*}
The $\xi_{c}(s)$ is an entire function of $s$ and satisfies the functional equation 
$\xi_{c}(s)=\xi_{c}(1-s).$\hfill\qed
\end{refproof}
%%%%%%%%%%%%%%%%%%%%%%%%%%%%%%%%%%%%%%%%%%%%%%%%%%%%
\begin{refproof}
This Corollary can be proved by using both Theorem 2 and the results of Biane, Pitman and Yor (To see, their table 1, P.442, 2001 ).
It follows from Biane, Pitman and Yor (2001) that the mixing density $f_z$ satisfies the relationship
$$
f_z(x)=x^{-\frac{3}{2}} \cdot f_z\left(\frac{1}{x}\right),\quad x>0
$$
and the moments of $f_z$ are finite for all orders and so are $f_{\hat{z}}$. In fact, the mixing density $f_{z}$ is equal to
$$
f_z(x)=\frac{1}{2 \sqrt{x}} H(\sqrt{x}), \quad x>0
$$
where the $H$ function is, for $y>0$,
$$
H(y)=4 y^2 \sum_{n=1}^{\infty}\left(2 \pi^2 n^4 y^2-3 \pi n^2\right) e^{-\pi n^2 y^2}.
$$
It follows from the Proposition 2.1 given by Biane, Pitman and Yor (2001) that a non-negative random variable $Y \geqslant 0$ has the density $y^{-1} H(y)$, for $y>0$. They obtained $Y^2 \stackrel{d}{=} \frac{\pi}{2} S_2$
( $\stackrel{d}{=}$ means the same distribution as), where the random variable $S_2 \geqslant 0$ has pdf $f_{S_2}$ and $f_{S_2}$ satisfies the relationship
$$
f_{S_2}(x)=\left(\frac{2}{\pi x}\right)^{\frac{5}{2}} \cdot f_{S_2}\left(\frac{4}{\pi^2 x}\right), \quad x>0 .
$$
Hence, it is easily to prove that
$$
f_z(x)=x^{-\frac{3}{2}} f_z\left(\frac{1}{x}\right), \quad x>0 .
$$
and the function $\xi_z(s)$ is equal to, for $s>1$,
$$
\begin{aligned}
\xi_z(s) & \equiv \int_0^{\infty} x^{\frac{s-1}{2}} f_z(x) d x \\
& =\int_0^{\infty} x^{\frac{s-1}{2}} \frac{1}{2 \sqrt{x}} H(\sqrt{x}) d x \\
& =\int_0^{\infty} y^{s-1} \frac{1}{2 y} H(y) \cdot 2 y d y \quad(\text { set } y=\sqrt{x}) \\
& =\int_0^{\infty} y^{s-1} H(y) d y=\int_0^{\infty} y^s \frac{1}{y} H(y) d y \\
& =E\left(Y^s\right) \\
& =2 \xi(s)
\end{aligned}
$$
the last equality follows from Biane, Pitman and Yor (2001, P.439). The other results of Corollary 3 follow immediately and the proof is complete.\hfill\qed
\end{refproof}
\begin{refproof}
First of all, we will prove that the function $f_w$ define on the $\mathbb{R}_{+}$is a probability density function, since this density function is a new one in the vein. The function
$$
f_w(x)=\frac{1}{2 \sqrt{x}} H_1(\sqrt{x}), \quad x>0
$$
is a pdf, which is equivalent to the function $H_1$ is a density function, since
$$
\begin{aligned}
\int_0^{\infty} f_w(x) d x & =\int_0^{\infty} \frac{1}{2 \sqrt{x}} H_1(\sqrt{x}) d x \quad\left(\text{Set}\;y=\sqrt{x}\right) \\
& =\int_0^{\infty} \frac{1}{2 y} H_1(y) 2 y d y \\
& =\int_0^{\infty} H_1(x) d x.
\end{aligned}
$$
We will prove that $H_1$ is a pdf as follows. Note that the function $H_1$ is equal to
$$
\begin{aligned}
H_1(x) & =\frac{d}{d x}\left(x^2 \frac{d}{d x} \theta_c(x)\right)\\
& =2 x \theta_c^{\prime}(x)+x^2 \cdot \theta_c^{\prime \prime}(x), \quad x>0
\end{aligned}
$$
Where $\theta_c(x)$ is a Jacobi's type theta function. Here, we choose the characteristic function $\hat{f}_c$ as given by Corollary 2, hence we have
$$
\begin{aligned}
\theta_c(x) & =\sum_{n=-\infty}^{\infty} \hat{f}_c(\sqrt{2 \pi} x n) \\
& =\sum_{n=-\infty}^{\infty} \frac{1}{\cosh (n \pi x)}, \quad x>0
\end{aligned}
$$
The Jacobi's type theta function $\theta_c$ satisfies the relationship
$$
\theta_c(x)=\frac{1}{x} \theta_c\left(\frac{1}{x}\right), \quad x>0
$$
and so is the function $H_1$, that is, the function $H_1$ satisfies the same functional equation as $\theta_c$ :
$$
H_1(x)=\frac{1}{x} H_1\left(\frac{1}{x}\right), \quad x>0
$$
(To see, Edwards, $1974\;\S\;10.3$, or Biane, Pitman and Yor 2001, P. 438 ).
Now, set $a=n \pi$, we have
$$
\begin{aligned}
H_1(x)= & 2 \sum_{n=1}^{\infty} a x \operatorname{sech}(a x)\left[a x\left(\tanh ^2(a x)-\operatorname{sech}^2(a x)\right)-2 \tanh (a x)\right] \\
= & 4 \sum_{n=1}^{\infty} \frac{a x\left[(a x-2) e^{2 a x}+(a x+2) e^{-2 a x}-6 a x\right]}{\left(e^{a x}+e^{-a x}\right)^3}
\end{aligned}
$$
It is clear that
$$
(a x-2) e^{2 a x}-6 a x>0 \text { for } n \geqslant 1 \text { and } x \geqslant 1.
$$
Hence $H_1(x)>0$ for all $x \geqslant 1$, and so $H_1(x)>0$ for all $x>0$, since
$$
H_1(x)=\frac{1}{x} H_1\left(\frac{1}{x}\right), \quad x>0 .
$$
Next the same reason as indicated by Riemann, the following integration by parts hold, for $s>1$
$$\begin{aligned} \int_0^{\infty} x^{s-1} H_1(x) d x  = & 2 \int_0^{\infty} x^{s-1} \frac{d}{d x}\left(x^2 \psi_c^{\prime}(x)\right) d x \\ = & 2 s(s-1) \int_0^{\infty} x^{s-1} \psi_c(x) d x \\ = & 2 \xi_c(s),\end{aligned}
$$
the last equality follows from Corollary 2, and it follows that
$$
\begin{aligned}
\int_0^{\infty} H_1(x) d x & =\lim _{s \rightarrow 1} \int_0^{\infty} x^{s-1} H_1(x) d x \\
& =\lim _{s \rightarrow 1} 2 \xi_C(s) \\
& =1 .
\end{aligned}
$$
Hence, we obtain that the function $H_1$ is a pdf and so is $f_w$. The mixing density $f_w$ also satisfies the relationship
$$
f_w(x)=x^{\frac{-3}{2}} \cdot f_w\left(\frac{1}{x}\right)
$$
This result follows from the relationship
$$
H_1(x)=\frac{1}{x} H_1\left(\frac{1}{x}\right).
$$
Similarly, we also have, for $s>1$
$$
\begin{aligned}
\xi_w(s) & =\int_0^{\infty} x^{\frac{s-1}{2}} f_w(x) d x \\
& =\int_0^{\infty} x^{s-1} H_1(x) d x \\
& =2 \xi_c(s)
\end{aligned}
$$
the other results follow immediately from Theorem 2, and the proof is complete.\hfill\qed
\end{refproof}
\begin{refproof}
This Corollary Can be proved by using both Theorem 2 and the result of Biane, Pitman and Your (To see, their Table 1, P.443,2001) as follows. They obtained that the Cosh-density $f_{C_1}$ on $\mathbb{R}_{+}$satisfies the relationship
$$
f_{C_1}(x)=\left(\frac{2}{\pi x}\right)^{\frac{3}{2}} f_{C_1}\left(\frac{4}{\pi^2 x}\right), \quad x>0
$$
and it is clearly that the following hold,
$$
T \stackrel{d}{=} \frac{\pi}{2} C_1
$$
hence it follows that the pdf $f_T$ satisfies
$$
f_T(x)=x^{-\frac{3}{2}} f_T\left(\frac{1}{x}\right), \quad x>0
$$
They also obtained
$$
\begin{aligned}
\xi_4(s) & =\int_0^{\infty} x^{\frac{s-1}{2}} f_T(x) d x \\
& =\left(\frac{4}{\pi}\right)^{\frac{s+1}{2}} \Gamma\left(\frac{s+1}{2}\right) L_{\chi_4}{(s)}
\end{aligned}
$$
where
$$
L_{\chi_4}(s)=\sum_{n=0}^{\infty} \frac{(-1)^n}{(2 n+1)^s}, \quad R_e(s)>0
$$
and $\xi_4$ is an entire function of $s \in \mathbb{C}$ and satisfies
$$
\xi_4(s)=\xi_4(1-s)
$$
and the other results follow immediately from Theorem 2. The proof is complete.\hfill\qed
\end{refproof}
\begin{refproof}
This Corollary follows immediately from Theorem 3. Since,
$$
g_2(x)=c \cdot x^{-\frac{3}{4}} h(x) h\left(\frac{1}{x}\right),\quad x>0
$$
Where $h(x)=e^{-a x^\alpha}$ is the Laplace-transform of a GGC and thus $g_2 \in B$ $($To see, Bondesson, 1992, P.75\;$)$.
It is clear that all the conditions of Theorem 3 hold, and the results follow. The proof is complete.\hfill\qed
\end{refproof}
\begin{refproof}
This corollary follows immediately from both Corollary 6 and the definition of generalized inverse Gaussian distribution, as follows its pdf is of the form
$$
\frac{\left(a/b\right)^{\frac{p}{2}}}{2 K_p(\sqrt{a b})} \cdot x^{p-1} \cdot e^{\frac{-(a x+\frac{b}{x})}{2}},\quad x>0
$$
Where $K_p$ is the modified Bessel function of the second kind, $a>0, b>0, p$ a real parameter. (To see, Jorgensen's lecture notes, 1982). Here, set $a=b>0$ and $p=\frac{1}{4}$, we get the mixing density
$$g_3(x)=\frac{1}{2 K_{\frac{1}{4}}(a)} x^{-\frac{3}{4}} \cdot e^{\frac{-a\left(x+\frac{1}{x}\right)}{2}},\quad x>0.$$
Note that the property of the generalized inverse Gaussian distribution implies
$$
\begin{aligned}
\xi_{G_3}(s) & =\int_0^{\infty} x^{\frac{s-1}{2}} \cdot g_3(x) d x \\
& =\int_0^{\infty} x^{\frac{s-1}{2}} \cdot \frac{1}{2 K_{\frac{1}{4}}(a)} x^{-\frac{3}{4}} \cdot e^{\frac{-a\left(x+\frac{1}{x}\right)} {2}} d x \\
& =\frac{K_{p_1}(a)}{K_{\frac{1}{4}}(a)} \cdot \int_0^{\infty} \frac{1}{2 K_{p_1}(a)} x^{p_1-1} \cdot e^{\frac{-a\left(x+\frac{1}{x}\right)}{2}} d x \\
& =\frac{K_{p_1}(a)}{K_{\frac{1}{4}}(a)} \cdot 1 \quad(\because p d f)
\end{aligned}
$$
Where $p_1=\frac{2 s-1}{4}$ and $a>0$.
The other results follow immediately from Corollary 6 and the proof is complete.\hfill\qed
\end{refproof}
\begin{refproof}
This Corollary follows immediately by using Corollary 6 and the Levy's formula (To see, Biane, Pitman and Yor, 2001, P. 444). The L\'{e}ry formula said
$$
\int_0^{\infty} \frac{a}{\sqrt{2 \pi t^3}} e^{-\frac{a^2}{2 t}} \cdot e^{-\lambda t} d t=e^{-a \sqrt{2 \lambda}}
$$
where $a>0$ and $\lambda>0$. Set $t=x^{-\frac{1}{2}}$ and $a^2=2 \lambda$, we get the mixing density $g_4$.
Now by using Corollary 6 with $\alpha=\frac{1}{2}$, the results follow immediately and the proof is complete.\hfill\qed
\end{refproof}
\end{document}